\documentstyle[12pt]{article}
\newfam\msbfam
\font\tenmsb=msbm10    \textfont\msbfam=\tenmsb
\font\sevenmsb=msbm7 \scriptfont\msbfam=\sevenmsb
\font\fivemsb=msbm5 \scriptscriptfont\msbfam=\fivemsb
\def\Bbb{\fam\msbfam \tenmsb}

\def\rr{{\Bbb R}}
\def\rz{{{\rr}^n}}
\def\zz{{\Bbb Z}}
\def\nn{{\Bbb N}}
\def\cc{{\Bbb C}}

\def\fz{\infty}
\def\az{\alpha}
\def\supp{{\rm{\ supp\ }}}
\def\loc{{\rm{\ loc\ }}}

\def\ez{\epsilon}
\def\bz{\beta}
\def\pz{\partial}
\def\tz{\theta}

\def\vz{\varphi}

\def\lz{\lambda}

\def\supp{{\rm supp}}
\def\loc{{\rm loc}}

\def\wz{\omega}
\def\l{\left}
\def\r{\right}

\def\dsum{\displaystyle\sum}

\def\dint{\displaystyle\int}
\def\dfrac{\displaystyle\frac}
\def\dsup{\displaystyle\sup}

\newtheorem{thm}{\hskip\parindent Theorem}

\newtheorem{lem}{\hskip\parindent Lemma}

\newtheorem{cor}{\hskip\parindent Corollary}

\begin{document}

\baselineskip=15pt
\renewcommand{\arraystretch}{2}
\arraycolsep=1.2pt

\title{A characterization of  weighted
local Hardy spaces } {\footnotetext{ \hspace{-0.65
cm} 2000 Mathematics Subject  Classification: 42B20, 42B25.\\
The  research was supported  by the NNSF (10971002) of China.\\}

\author{ Lin Tang }
\date{}
\maketitle

{\bf Abstract}\quad  In this paper, we  give a  characterization
of weighted local Hardy spaces $h^1_\wz(\rz)$ associated with
local weights  by using the truncated Reisz transforms, which
generalizes the corresponding result of Bui in \cite{b}.

\bigskip

\begin{center}{\bf 1. Introduction }\end{center}
The theory of local Hardy space plays an important role in various
fields of analysis and partial differential equations; see
\cite{g,b}.  Bui \cite{b} studied the weighted version $h_w^p$ of
the local Hardy space $h^p$ considered by Goldberg \cite{g}, where
the weight $\wz$ is assumed to satisfy the condition $(A_\fz)$ of
Muckenhoupt.  R. Vyacheslav \cite{v} introduced and studied some
properties of the weighted local Hardy space $h_\wz^p$ spaces
with weights that are locally in $A_p$ but may grow or decrease
exponentially. Recently, the author \cite{t} established the
weighted atomic decomposition characterizations of weighted local
Hardy space $h_\wz^p$ with local weights.

The main purpose of this paper is to give a  characterization of
weighted local Hardy spaces $h^1_\wz(\rz)$ associated with local
weights by using the truncated Reisz transforms.

Throughout this paper, $C$ denotes the constants that are
independent of the main parameters involved but whose value may
differ from line to line. Denote by $\nn$ the set $\{1,2,\cdots\}$
and by $\nn_0$ the set $\nn\cup\{0\}$. By $A\sim B$, we mean that
there exists a constant $C>1$ such that $1/C\le A/B\le C$.

\begin{center}{\bf 2.  Statement of the main result }\end{center}
We first introduce  weight classes $A_p^{loc}$ from \cite{v}.

Let $Q$ run through all cubes in $\rz$ (here and below only cubes
with sides parallel to the coordinate axes are considered), and
let $|Q|$ denote the volume of $Q$. We define the weight class
$A_p^{loc}(1<p<\fz$) to consists of all nonnegative locally
integral functions $\wz$ on $\rz$ for which
$$A_p^{loc}(\wz)=\dsup_{|Q|\le 1}\dfrac
1{|Q|^p}\dint_Q\wz(x)dx\l(\dint_Q\wz^{-p'/p}(x)dx\r)^{p/p'}<\fz ,
\  1/p+1/p'=1.\eqno(2.1)$$ The function $\wz$ is said to belong to
the weight class of $A_1^{loc}$ on $\rz$ for which
$$A_1^{loc}(\wz)=\dsup_{|Q|\le 1}\dfrac
1{|Q|}\dint_Q\wz(x)dx\l(\dsup_{y\in
Q}[\wz(y)]^{-1}\r)<\fz.\eqno(2.2)$$ {\bf Remark}: For any $C>0$ we
could have replaced $|Q|\le 1$ by $|Q|\le C$ in (2.1) and (2.2).

 In what follows,
$Q(x,t)$ denotes the cube centered at $x$ and of the sidelength
$t$. Similarly, given $Q=Q(x,t)$ and $\lz>0$, we will write $\lz
Q$ for the $\lz$-dilate cube, which is the cube with the same
center $x$ and with sidelength $\lz t$. Given a Lebesgue
measurable set $E$ and a weight $\wz$,  let $\wz(E)=\int_E\wz dx$.
For any $\wz\in A_\fz^{loc}$, $L^p_\wz$ with $p\in(0,\fz)$ denotes
the set of all measurable functions $f$ such that
$$\|f\|_{L^p_\wz}\equiv \l(\dint_\rz|f(x)|^p\wz(x)dx\r)^{1/p}<\fz,$$
and $L^\fz_\wz=L^\fz$. The space $L_\wz^{1,\fz}$ denotes the set
of all measurable function $f$ such that
$$\|f\|_{L_\wz^{1,\fz}}\equiv\dsup_{\lz>0}\lz\cdot \wz(\{x\in\rz:
|f(x)|>\lz\})<\fz.$$ We define the local Hardy-Littlewood maximal
operator by
$$M^{loc}f(x)=\dsup_{x\in Q:|Q|<1}\dfrac 1{|Q|}\dint_Q|f(y)|dy.$$

 Similar to the
classical $A_p$ Muckenhoupt weights, we give some properties for
weights $\wz\in A^\loc_\fz:=\bigcup_{1\le p<\fz} A^{loc}_p$.
\begin{lem}\label{l2.1.}\hspace{-0.1cm}{\rm\bf 2.1.}\quad
Let $1\le p<\fz$, $\wz\in A_p^{loc}$, and $Q$ be a  unit cube,
i.e. $|Q|=1$. Then there exists a  $\bar\wz\in A_p$ so that
$\bar\wz=\wz$ on $Q$ and
\begin{enumerate}
\item[(i)]$A_p(\bar \wz)\le CA_p^{loc}(\wz).$ \item[(ii)]if
$\wz\in A_p^{\loc}$, then there exists $\ez>0$ such that $\wz\in
A_{p-\ez}^{loc}(\wz)$ for $p>1$. \item[(iii)]If $ 1\le
p_1<p_2<\fz$, then $A_{p_1}^{loc}\subset A_{p_2}^{loc}$.
\item[(iv)] $\wz\in A_p^{loc}$ if and only if $\wz^{-\frac
1{p-1}}\in A_{p'}^{loc}$. \item[(v)] If $\wz\in A_p^{loc}$ for
$1\le p<\fz$, then
$$\wz(tQ)\le exp(c_\wz t)\wz(Q)\quad (t\ge 1, |Q|=1).$$
\item[(vi)] the local Hardy-Littlewood maximal operator $M^{loc}$
is bounded on $L^p_\wz$ if $\wz\in A_p^{loc}$ with $p\in (1,\fz)$.
\item[(vii)]$M^{loc}$ is bounded from $L^1_\wz$ to $L^{1,\fz}_\wz$
if $\wz\in A_1^{loc}$.
\end{enumerate}
\end{lem}

We remark that Lemma is also true for $|Q|>1$ with $c$ depending
now on the size of $Q$. In addition, it is easy to see that
$A_p(Munckenhoupt\ weight)\subset A_p^{loc}$ for $p\ge 1$ and
$e^{c|x|},\ (1+|x|\ln^\az(2+|x|))^\bz\in A_1^{\loc}$ with $\az\ge
0,\bz\in\rr$ and $c\in\rr$.

Let ${\cal N}$ denote the class of $C^\fz$-functions $\vz$ on
$\rz$, supported on the cube $Q(0,1)$ of center zero and half-side
one whose mean value is not equal to zero. For
$t>0$, let $\vz_t=t^{-2n}\vz(z/t)$.

Given a distribution $f$, let $\vz\in{\cal N}$, define the smooth maximal function by
$${\cal M}f(z)=\dsup_{0<t<1}|\vz_t*f(z)|.$$
Follows from \cite{t}, we introduce the following weighted atoms.

Let $\wz\in A_1^{loc}$. A function $a$ on $\rz$ is said to be a
$(1,q)_\wz$-atom for $1<q\le\fz$ if
\begin{enumerate}
\item[(i)] $\supp \ a\subset Q$, \item[(ii)]
$\|a\|_{L^q_\wz(\rz)}\le [\wz(Q)]^{1/q-1}$. \item[(iii)]$\dint_\rz
a(x) dx=0$  if $|Q|<1$.
\end{enumerate}
Moreover, we call $a$ is a $(1,q)_\wz$ single atom if
$\|a\|_{L^q_\wz(\rz)}\le [\wz(\rz)]^{1/q-1}$. we introduce weighted
local Hardy spaces via smooth maximal functions and weighted local
Hardy spaces. Moreover, we study some properties of these spaces.

The weighted local Hardy space is defined by
$$H_{\wz}^1(\rz)\equiv\l\{f\in {\cal D}'(\rz): {\cal M}(f)\in
L_\wz^1(\rz)\r\}.$$ Moreover, we define
$\|f\|_{h^1_{\wz}(\rz)}\equiv\|{\cal
M}(f)\|_{L^1_\wz(\rz)}$. In \cite{t}, the author proved that
\begin{thm}\label{tA.}\hspace{-0.1cm}{\rm\bf A.}\quad
Let $\wz\in A_1^{loc}$ and $1<q\le \fz$, then for any $f\in h^1_{\wz}(\rz)$, there
exists numbers $\lz_0$ and $\{\lz_i^k\}_{k\in\zz,i}\subset \cc$,
$(1,q)_\wz$-atoms $\{a_i^k\}_{k\in\zz,i}$ with radius $r\le2$
and single atom $a_0$ such that
$$f=\dsum_{k\in\zz}\dsum_{i}\lz_i^ka_i^k+\lz_0a_0,$$ where the series converges
almost everywhere and in ${\cal D}'(\rz)$, moreover, there exists
a positive constant $C$, independent of $f$, such that
$\dsum_{k\in\zz,i}|\lz_i^k|^p+|\lz_0|^p\le
C\|f\|_{h^1_{\wz}(\rz)}$.
\end{thm}

Let $\Phi$ be a non-negative, radial and $C^\fz$-function on $\rz$ with compact support $B(0,2)$ and
 $\Phi\equiv 1$ on $B(0,1)$. Define  the truncated Reisz
transforms by
$$R_jf(x)=\dint_\rz K_j(x-y)f(y)dy,\ K_j(z)=\dfrac {z_j}{|z|^{n+1}}\Phi(z),\ j=1,\cdots,n.$$
Now let us state the main result of this paper.
\begin{thm}\label{t1.}\hspace{-0.1cm}{\rm\bf 1.}\quad
Let $\wz\in A_1^{loc}$. Then a function $f$ is in $h_\wz^1(\rz)$
if and only if $f\in L_\wz^1(\rz)$ and $R_jf\in L^1_\wz(\rz)$,
$j=1,\cdots, n$. More precisely,
$$\|f\|_{h_\wz^1(\rz)}\sim \|f\|_{L^1_\wz(\rz)}+\dsum_{j=1}^n\|R_j
f\|_{L^1_\wz(\rz)}.$$
\end{thm}
We remark that if $\wz\in A_1$, then Theorem 1 has been proved in
\cite{b}, that is,

\begin{thm}\label{tB.}\hspace{-0.1cm}{\rm\bf B.}\quad
Let $\wz\in A_1$. Then a function $f$ is in $h_\wz^1(\rz)$ if and
only if $f\in L_\wz^1(\rz)$ and $R_jf\in L^1_\wz(\rz)$,
$j=1,\cdots, n$. More precisely,
$$\|f\|_{h_\wz^1(\rz)}\sim \|f\|_{L^1_\wz(\rz)}+\dsum_{j=1}^n\|R_j
f\|_{L^1_\wz(\rz)}.$$
\end{thm}

\begin{center} {3. Proof of Theorem 1 }\end{center}
Theorem 1 will be deduced by  the following lemmas.
\begin{lem}\label{l3.1.}\hspace{-0.1cm}{\rm\bf 3.1.}\quad
Let $\wz\in A_1^{loc}$. Then $$\|f\|_{h_\wz^1(\rz)}\le C(
\|f\|_{L^1_\wz(\rz)}+\dsum_{j=1}^n\|R_j
f\|_{L^1_\wz(\rz)}).\eqno(3.1)$$
\end{lem}
Proof.\quad We will borrow some idea from \cite{mpr}. Let $Q$ is
an unit cube, $\chi'_{3Q}$ is a $C^\fz_0$ nonnegative function
supported in $4Q$ and $\chi'_{3Q}=1$ on $3Q$. By Lemma 2.1, we can
set $\bar\wz\in A_p$ so that $\bar\wz=\wz$ on $14Q$. Fix $\vz\in{\cal N}$, by Theorem B,
we have
$$\begin{array}{cl}
\|\dsup_{0<t<1}|\vz_t*f\|_{L^1_\wz(Q)}&=
\|\dsup_{0<t<1}|\vz_t*(f\chi'_{3Q})\|_{L^1_{\bar \wz}(\rz)}\\
&\le C\|f\chi'_{3Q}\|_{h^1_{\bar \wz}(\rz)}\\
&\le C\l(\|f\chi'_{3q}\|_{L^1_{\bar \wz}(\rz)}+\dsum_{j=1}^n\|R_j(
f\chi'_{3Q})\|_{L^1_{\bar \wz}(\rz)}\r).
\end{array}\eqno(3.2)$$
On the other hand, by the properties of $A_1^{loc}$, we obtain
$$\begin{array}{cl}
\|R_j( f\chi'_{3Q})&-\chi'_{3Q}R_j( f)\|_{L^1_{\bar\wz}(\rz)}\\&\le
\|\dint|R_j(z-y)[\chi'_{3Q}(y)-\chi'_{3Q}(z)]f(y)|\chi'_{12I}(y)dy\|_{L^1_{\bar\wz}(\rz)}\\
&\le C\dint_\rz\bar\wz(z)\dint_\rz|R_j(z-y)||z-y||f(y)|\chi'_{12I}(y)dydz\\
&\le C\|f\|_{L^1_{\bar\wz}(14Q)}.
\end{array}\eqno(3.3)$$
Combing (3.2) and (3.3), we obtain
$$
\|\dsup_{0<t<1}|\vz_t*f|\|_{L^1_\wz(Q)}\le C\l(\|f\|_{L^1_{
\wz}(14Q)}+ \dsum_{j=1}^n\|R_j( f)\|_{L^1_{ \wz}(6Q)})\r).$$
Summing on $Q$, we obtain (3.1).

\begin{lem}\label{l3.2.}\hspace{-0.1cm}{\rm\bf 3.2.}\quad
Let $R_j$ be as above, then
\begin{enumerate}
\item[(i)]$\|R_j f\|_{L^p_\wz(\rz)}\le
C_{p,\wz}\|f\|_{L^p_\wz(\rz)}$ for $1<p<\fz$ and $\wz\in
A_p^{loc}$. \item[(ii)]$\|R_j f\|_{L^{1,\fz}_\wz(\rz)}\le
C_\wz\|f\|_{L^1_\wz(\rz)}$ for $\wz\in A_1^{loc}$.
\end{enumerate}
\end{lem}
Proof.  We first note that for $\wz\in A_p$ the inequality (i) is
known to be true, see \cite{gr}. For $\wz\in A_p^{loc}$, by Lemma
2.1 (i) for any unit cube $Q$ there is a $\bar\wz\in A_p$ so that
$\bar\wz=\wz$ on $6Q$. Then
$$\begin{array}{cl}
\|R_jf\|_{L^p_\wz(Q)}&= \|R_j(\chi_{6Q}f)\|_{L^p_\wz(Q)}\\
&\le \|R_j(\chi_{6Q}f)\|_{L^p_{\bar\wz}(Q)}\\
&\le C\|(\chi_{6Q}f)\|_{L^p_{\bar\wz}(\rz)}\\
&\le C\|f\|_{L^p_{\bar\wz}(6Q)}.
\end{array}$$
Summing over all dyadic unit $Q$ gives (i).

For (ii), similar to (i), note that for $\wz\in A_1$ the
inequality (ii) is known to be true, see \cite{c}. Since $\wz\in
A_p^{loc}$, by Lemma 2.1 (i) for any unit cube $Q$ there is a
$\bar\wz\in A_1$ so that $\bar\wz=\wz$ on $6Q$. Then for any
$\lz>0$
$$\begin{array}{cl}
\wz(\{x\in Q: |R_jf(x)|>\lz\})&\le \wz(\{x\in Q: |R_j(\chi_{6Q}f)(x)|>\lz\})\\
&= \bar\wz(\{x\in Q: |R_j(\chi_{6Q}f)(x)|>\lz\})\\
&\le C\lz^{-1}\|(\chi_{6Q}f)\|_{L^1_{\bar\wz}(\rz)}\\
&= C\lz^{-1}\|f\|_{L^1_{\wz}(6Q)}.
\end{array}$$
Summing over all dyadic unit $Q$ gives (ii).
\begin{lem}\label{l3.3.}\hspace{-0.1cm}{\rm\bf 3.3.}\quad
Let $\wz\in A_1^{loc}$. Then $$\|R_jf\|_{h_\wz^1(\rz)}\le C
\|f\|_{h^1_\wz(\rz)}.\eqno(3.4)$$
\end{lem}
Proof:\quad  We first fix a function $\vz\in{\cal N}$. Let $a(x)$ be a $(1,2)$
atom in $h^1_\wz(\rz)$, supported in a cube $Q$ centered at $y_0$
and sidelength $r\le2$, or $a(x)$ is a $(1,2)$ single atom. To prove the
(iii), by Theorem A and Theorem 6.2 in \cite{t}, it is enough to show that
$$\|{\cal M}(R_ja)\|_{L^1_\wz(\rz)}\le C,\eqno(3.5)$$
where $C$ is independent of $a$.

If $a$ is a single atom, by $L^2_\wz(\rz)$ boundedness of ${\cal
M}$ and $R_j$, then
$$\|{\cal M}(R_j a)\|_{L^1_\wz(\rz)}\le
C\|R_ja\|_{L^2_\wz(\rz)}\wz(\rz)^{1/2}\le C.$$

Next we always assume that $a$ is an atom in $h^1_\wz(\rz)$,
supported in a cube $Q$ centered at $y_0$ and sidelength $r\le2$.

We first consider the atom $a$ with sidelength $1\le r\le 2$. Then
by $L^2_\wz(\rz)$ of the operators ${\cal M}$ and $R_j$(see
Lemma 3.2), we have
$$\begin{array}{cl}
\|{\cal M}(R_j a)\|_{L^1_\wz(\rz)}&= \dint_{8Q} {\cal  M}(R_j
a)(x)\wz(y)dy\\
&\le C\wz(8Q)^{1/2}\|a\|_{L^2_\wz(\rz)}\le C.
\end{array}$$
If $r<1$, we write
$$\begin{array}{cl}
\|{\cal M}(R_j a)\|_{L^1_\wz(\rz)}&= \dint_{2Q} {\cal M}(R_j
a)(x)\wz(y)dy+\dint_{\rz\setminus 2Q} {\cal M}(R_j
a)(x)\wz(y)dy\\
&:=I+II.
\end{array}$$
For $I$, by $L^2_\wz(\rz)$ boundedness of the operators ${\cal M}$ and
$R_j$, we have
$$I\le \wz(2Q)^{1/2}\|a\|_{L^2_\wz(\rz)}\le C.$$
We now estimate $II$. Let $x\notin 2Q$. For $t>0$ we define the smooth functions
$$R_j^t=\vz_t*K_j$$ and we observe that they satisfy
$$\dsup_{0<t<1}|\pz^\bz K_j^t(x)|\le C|x-y_0|^{-n-|\bz|}\chi_{\{|x-y_0|\le 8n\}}(x)\eqno(3.6)$$
for all $|\bz|\le 1$; see their proof in page 507 of \cite{gra}.

Now note that if $x\not\in 2Q$  and  and $y\in Q$, then $|x-y_0|\ge 2|y-y_0|$ stays away
from $y_0$ and $K_j(x-y)$ is well defined. We have
$$R_ja*\vz_t(x)=(a*K_j^t)(x)=\dint_Q K^t_j(x-y)a(y)dy.$$
Using the cancellation of atoms we deduce
$$\begin{array}{cl}
R_ja*\vz_t(x)&=\dint_Q K^t_j(x-y)a(y)dy\\
&=\dint_Q \l[K^t_j(x-y)-K^t_j(x-y_0)\r]a(y)dy\\
&=\dint_Q \l[\dsum_{|\bz|=1}(\pz^\bz K_j^t(x-y_0-\tz_y(y-y_0))y^\bz\r]a(y)dy
\end{array}$$
for some $0\le\tz_y\le 1$. Using that $|x-y_0|\ge 2|y-y_0|$ and (3.6) we get
$$\begin{array}{cl}
R_ja*\vz_t(x)&\le C|x-y_0|^{-n-1}\chi_{\{|x-y_0|\le 8n\}}(x)\dint_Q|a(y)||y|dy\\
&\le C\dfrac {r^{n+1}}{|x-y_0|^{n+1}}\wz(Q)^{-1}\chi_{\{|x-y_0|\le 8n\}}(x).
\end{array}\eqno(3.7)$$
By (3.7) and using properties of $A_1^{loc}$, we obtain
$$\begin{array}{cl}II&\le C\dint_{2r\le|x-y_0|\le 8n}\dfrac
{r^{n+1}}{|x-y_0|^{n+1}}\wz(Q)^{-1}\wz(x)dx\\
&\le
C\dfrac{|Q|}{\wz(Q)}\dsum_{k=1}^{k_0}2^{-k}\dfrac{\wz(2^kQ)}{|2^kQ|}\le
C,
\end{array}$$ where $k_0$ is an integer such that $8n\le 2^{k_0}\le
16n$.

 Thus, (3.5)  holds. Hence, the proof is complete.

 Next, we study weighted $h^1_\wz(\rz)$ boundedness for strongly
singular integrals.

Given a  real number $\tz>0$ and a smooth radial cut-off function
$v(x)$ supported in the ball $\{x\in\rz: \ |x|\le 2\}$, we
consider the strongly singular kernel
$$k(x)=\dfrac{e^{i|x|^{-\tz}}}{|x|^n}v(x).$$
Let us denote by $Tf$ the corresponding strongly singular integral
operator:
$$Tf(x)=p.v\dint_\rz k(x-y)f(y)dy.$$
This operator has been studied by several authors, see \cite{h},
\cite{w} and \cite{ghst}. In particular, S. Chanillo \cite{c}
established  the weighted $L^p_\wz(\rz)$($\wz\in A_p,1<p<\fz$ )
and $H^1_\wz(\rz)$($\wz\in A_1$) boundedness for strongly singular
integrals. The author \cite{t} proved the following results for
the strongly singular integrals.
\begin{thm}\label{tC.}\hspace{-0.1cm}{\rm\bf C.}\quad
Let $T$ be strongly singular integral operators, then
\begin{enumerate}
\item[(i)]$\|Tf\|_{L^p_\wz(\rz)}\le C_{p,\wz}\|f\|_{L^p_\wz(\rz)}$
for $1<p<\fz$ and $\wz\in A_p^{loc}$.
\item[(ii)]$\|Tf\|_{L^{1,\fz}_\wz(\rz)}\le
C_\wz\|f\|_{L^1_\wz(\rz)}$ for $\wz\in A_1^{loc}$.
\item[(iii)]$\|Tf\|_{L^1_\wz(\rz)}\le C_\wz\|f\|_{h^1_\wz(\rz)}$
for $\wz\in A_1^{loc}$.
\end{enumerate}
\end{thm}

 Theorem 1, Lemma 3.3 and (iii) in Theorem C imply immediately that
\begin{cor}\label{c1.}\hspace{-0.1cm}{\rm\bf 1.}\quad
Let $T$ be strongly singular integral operators, then
$$\|Tf\|_{h^1_\wz(\rz)}\le C_\wz\|f\|_{h^1_\wz(\rz)}$$ for $\wz\in
A_1^{loc}$.
\end{cor}

 LMAM, School of Mathematical  Science

 Peking University

 Beijing, 100871

 P. R. China

\bigskip

 E-mail address:  tanglin@math.pku.edu.cn


\begin{center} {\bf References}\end{center}
\begin{enumerate}
 \vspace{-0.3cm}
\bibitem[1]{b} H. Bui,
Weighted Hardy spaces, Math. Nachr. 103 (1981), 45--62.
\vspace{-0.3cm}
\bibitem[2]{c} S. Chanillo,
Weighted norm inequalities for strongly singular convolution
operators, Trans. Amer. Math. Soc. 281(1984), 77-107.
 \vspace{-0.3cm}
\bibitem[3]{gc} J. Garc\'ia-Cuerva,
Weighted $H^p$ spaces, Dissertationes Math. 162(1979), 63.
\vspace{-0.3cm}
\bibitem[4]{ghst} J. Garc\'ia-Cuerva, E. Harboure, S. Segovia and J. L. Torrea,
Weighted norm inequalities for commutators of strongly singular
integral, Indiana. Univ. Math. J. 40(1991), 1397-1420.
 \vspace{-0.3cm}
\bibitem[5]{gr} J. Garc\'ia-Cuerva and J. Rubio de Francia,
Weighted norm inequalities and related topics, Amsterdam- New
York, North-Holland, 1985. \vspace{-0.3cm}
\bibitem[6]{g} D. Goldberg, A local version of real Hardy spaces,
Duke Math. 46(1979), 27-42.
\vspace{-0.3cm}
\bibitem[7]{gra} L. Grafakos,
Classical and modern fourier analysis, 2004.
\vspace{-0.3cm}
\bibitem[8]{h}I. Hirschman,
Multiplier transformations, Duke Math. J. 26(1959), 222-242.
\vspace{-0.3cm}
\bibitem[9]{mpr}G. Mauceri, M. Picardello and F. Ricci,
A Hardy space associated with twisted convolution, Adv. Math.
39(1981), 270-288.
 \vspace{-0.3cm}
\bibitem[10]{t}L. Tang,
Weighted local Hardy spaces and their applications, preprint.
\vspace{-0.3cm}
\bibitem[11]{v} R. Vyacheslav, Littlewood-Paley
theory and function spaces with $A^{\rm loc}_p$ weights, Math.
Nachr. 224 (2001), 145--180. \vspace{-0.3cm}
\bibitem[12]{w}  S. Wainger,
Special trigonometric series in $k$ dimensions, Mem. Amer. Math.
Soc. 59(1965).
\end{enumerate}
\end{document}